\documentclass[fleqn]{mat01}
\usepackage{times,mathtimy,amssymb,latexsym,amsbsy}
\begin{document}

\setcounter{page}{103}
\firstpage{103}

\newtheorem{theor}{\bf Theorem}
\newtheorem{pot}{\it Proof of Theorem}
\newtheorem{propo}{\rm PROPOSITION}

\def\theorr{\trivlist \item[\hskip \labelsep{\bf Theorem S.}]}
\def\theorre{\trivlist \item[\hskip \labelsep{\bf Theorem.}]}

\title{Isometric multipliers of $\pmb{L^p (G,X)}$}

\markboth{U~B~Tewari and P~K~Chaurasia}{Isometric multipliers of ${L^p (G,X)}$}

\author{U~B~TEWARI and P~K~CHAURASIA}

\address{Department of Mathematics, Indian Institute of Technology,
Kanpur~208~016, India\\
\noindent E-mail: ubtewari@iitk.ac.in; praveenc@iitk.ac.in}

\volume{115}

\mon{February}

\parts{1}

\pubyear{2005}

\Date{MS received 21 January 2004}

\begin{abstract}
Let $G$ be a locally compact group with a fixed right Haar measure and
$X$ a separable Banach space. Let $L^p(G,X)$ be the space of $X$-valued
measurable functions whose norm-functions are in the usual $L^{p}$. A
left multiplier of $L^p(G,X)$ is a bounded linear operator on $
L^p(G,X)$ which commutes with all left translations. We use the
characterization of isometries of $L^p(G,X)$ onto itself to characterize
the isometric, invertible, left multipliers of $L^p(G,X)$ for $1\leq p
<\infty$, $p\neq 2$, under the assumption that $X$ is not the
$\ell^p$-direct sum of two non-zero subspaces. In fact we prove that if
$T$ is an isometric left multiplier of $L^p(G,X)$ onto itself then there
exists a $y \in G$ and an isometry $U$ of $X$ onto itself such that
$Tf(x)= U(R_yf)(x)$. As an application, we determine the isometric left
multipliers of $L^1 \bigcap L^p (G,X)$ and $L^1 \bigcap C_0 (G,X)$ where
$G$ is non-compact and $X$ is not the $\ell^p$-direct sum of two
non-zero subspaces. If $G$ is a locally compact abelian group and $H$ is
a separable Hilbert space, we define $A^p (G,H) = \{f\in
L^1(G,H)\hbox{:}\ \hat{f}\in L^p(\Gamma,H)\}$ where $\Gamma$ is the dual
group of $G$. We characterize the isometric, invertible, left multipliers
of $A^p (G,H)$, provided $G$ is non-compact. Finally, we use the
characterization of isometries of $C(G,X)$ for $G$ compact to determine
the isometric left multipliers of $C(G,X)$ provided $X^*$ is strictly
convex.
\end{abstract}


\maketitle

\section{Introduction}

Let $G$ be a locally compact group with right Haar measure $\mu$.
Suppose $X$ is a separable Banach space. If $1\leq p<\infty$, let
$L^p(G,X)$ be the space of $X$-valued measurable functions $F$ such that
$\int_G \|F(x)\|^p\ {\rm d}\mu < \infty$. The $p$-norm of $F$ is defined
by $[\int_G\|F\|^p\ {\rm d}\mu]^{1/p}$. In case $X$ is a
one-dimensional complex Banach space, $L^p(G,X)$ is denoted by $L^p(G)$.

The left and right translation operators $L_g$ and $R_g$ are defined by
$(L_g F)(x)= F(gx)$ and $(R_g F)(x) = F(xg)$. A left multiplier of
$L^p(G,X)$ is a bounded linear operator on $L^p(G,X)$ which commutes
with all left translations. The main result of this paper gives a
characterization of the isometric, invertible, left multipliers of
$L^p(G,X)$ for $1\leq p <\infty, p\neq 2$, under the assumption that $X$
is not the $\ell^p$-direct sum of two non-zero subspaces. More precisely
we shall prove the following theorem.

\begin{theor}[\!] Let $G$ be a locally compact group and $T$ an
isometric{\rm ,} invertible{\rm ,} left multiplier on $L^p(G,X)$ for
$1\leq p<\infty, p\neq 2$. Suppose that $X$ is not the $\ell^p$-direct
sum of two non-zero subspaces. Then there exists an isometry $U$ of $X$
onto itself and $y\in G$ such that $T$ is of the form
\begin{equation*}
(TF)(x) = U R_y F(x).
\end{equation*}
\end{theor}
Wendel \cite{7} proved this result for $L^1(G)$ in 1952. Later
Strichartz \cite{6} and Parrot \cite{4} proved it for $L^p (G)$ if
$1\leq p<\infty, p\neq 2$.

Let $G$ be a non-compact locally compact group. If $f\in L^1 \bigcap L^p
(G,X)$, we define $\|f\| = \|f\|_1 + \|f\|_p$. Then $L^1 \bigcap L^p
(G,X)$ is a Banach space with this norm. Similarly for $f\in L^1 \bigcap
C_0 (G,X)$, we define $\|f\| = \|f\|_1 + \|f\|_\infty$. Then $L^1 \bigcap
C_0 (G,X)$ is a Banach space. In both cases, we shall show that if $T$
is an isometric, invertible, left multiplier, then $T$ is of the form
\begin{equation*}
(Tf)(x) = U R_y f(x).
\end{equation*}
If $G$ is a locally compact abelian group and $H$ is a separable Hilbert
space, we define $A^p (G,H) = \{f\in L^1(G,H)\hbox{:}\ \hat{f}\in
L^p(\Gamma,H)\}$ where $\Gamma$ is the dual group of $G$. For $f\in A^p
(G,H)$, we define $\|f\| = \|f\|_1 + \|\hat{f}\|_p$. $A^p (G,H)$ is a
Banach space with this norm. We will prove that if $T$ is an isometric,
invertible, left multiplier, then $T$ is of the form
\begin{equation*}
(Tf)(x) = U R_y f(x).
\end{equation*}
Let $G$ be a compact group and $X$ be a separable Banach space. $C(G,X)$
denotes the Banach space of continuous $X$-valued functions. Using the
characterization of isometries of $C(G,X)$, we will prove that if $T$ is
an isometric, invertible, left multiplier, then $T$ is of the form
\begin{equation*}
(TF)(x) = U R_y F(x).
\end{equation*}
provided $X^*$ is strictly convex.

\section{Preliminaries}

Let $(\Omega, \Sigma, \mu)$ be a measure space. Suppose $\Sigma^\prime$
is the $\sigma$-ring generated by the sets of $\sigma$-finite measure. A
mapping $\Phi$ of $\Sigma^\prime$ onto itself, defined modulo null sets,
is said to be a regular set isomorphism if
\begin{enumerate}
\renewcommand\labelenumi{\arabic{enumi}.}
\leftskip -.2pc
\item $\Phi (A\backslash A^\prime) = \Phi (A) \backslash \Phi
(A^\prime)$ for $A, A^\prime \in \Sigma^\prime$.\vspace{.3pc}

\item $\Phi(\bigcup_{n=1}^\infty A_n) = (\bigcup_{n=1}^\infty \Phi
(A_n)$, where $\{A_n\}$ is a sequence of disjoint sets in $\Sigma^\prime$.\vspace{.3pc}

\item $\mu (\Phi (A)) = 0$ iff $\mu (A) =0$.
\end{enumerate}
\looseness -1 A regular set isomorphism induces a linear map on $X$-valued measurable
functions. If $A\in \Sigma^\prime$ and $x\in X$, define $\Phi (\chi_{A})(x)
= \chi_{\Phi(A)}x$ where $\chi_A$ is the characterstic function of $A$.
This extends linearly to simple functions. Let $f$ be an $X$-valued
measurable function. Then there exists a sequence $\{f_n\}$ of simple
functions converging to $f$ in measure. Then $\{\Phi(f_n \}$ is a Cauchy
sequence in measure and hence converges to a measurable function
$\Phi(f)$. It is easy to show that $\Phi(f)$ depends only on $f$ and not
on the particular sequence $\{f_n\}$.

We also note that any $\Sigma^\prime$-measurable function is also
$\Sigma$-measurable and any $\Sigma$-measurable function with
$\sigma$-finite support is $\Sigma^\prime$-measurable. Thus the spaces of
$\Sigma^\prime$ and $\Sigma$ measurable functions with $\sigma$-finite
support coincide.

If $\Phi$ is a regular set isomorphism, define a measure $\nu$ by
$\nu(A) = \mu(\Phi^{-1} (A))$. The measure $\nu$ is absolutely continuous
with respect to $\mu$. Let $h = (({\rm d}\nu)/{\rm d}\mu)^{1/p}$. It is
easy to see that $h$ is a function on $\Omega$ whose restriction to any
measurable set of $\sigma$-finite measure is measurable. Further, if
$f\in L^p(\Omega, X)$, then $h\Phi(f) \in L^p(\Omega,X)$ and
$\|h\Phi(f)\|_p = \|f\|_p$.

We say that a Banach space $X$ is the $\ell^p$-direct sum of two Banach
spaces $X_1$ and $X_2$ if $X$ is isometrically isomorphic to
$X_1\bigoplus X_2$ where the norm on the direct sum is given by
$\|x_1\bigoplus x_2\| = \{\| x_1\|^p + \|x_2\|^p\}^{1/p}$.

Our main tool for the proof of the main result is a theorem of Sourour
\cite{5}. We state it in a form slightly different from that of
\cite{5}, but virtually no modification of the proof given there is
necessary. The assumption that $\Omega$ is $\sigma$-finite is not needed
for our conclusion because every function in $L^p(\Omega,X)$ has
$\sigma$-finite support.

\begin{theorr}{\it Let $(\Omega,\Sigma,\mu)$ be a measure space and $T$
be an isometry of $L^p(\Omega,X)$ onto itself. Suppose $X$ is not the
$\ell^p$-direct sum of two non-zero Banach spaces. Then there exists a
regular set isomorphism $\Phi$ of $\Sigma^\prime$ onto itself{\rm ,} a
measurable function $h$ on $\Omega$ and a strongly measurable map $S$ of
$\Omega$ into the Banach space of bounded linear maps of $X$ into $X$
with $S(t)$ a surjective isometry of $X$ for almost all $t\in
\Omega${\rm ,} such that
\begin{equation*}
T(F(t)) = S(t) h(t) \Phi(F)(t)
\end{equation*}
for $F\in L^p(\Omega,X)$ and almost all $t\in \Omega$.}
\end{theorr}

\section{Isometric multipliers of $\pmb{L^p(G,X)}$}

In this section we characterize the isometric, invertible, left
multipliers of $L^p(G,X)$.

\begin{pot}{\rm Let $T$ be an isometric, invertible, left multiplier of
$ L^p(G,X)$. It follows from Theorem~S that
\begin{equation*}
TF(t) = h(t) S(t)\Phi(F)(t)\quad\ \hbox{a.e.}
\end{equation*}
for every $F\in L^p(G,X)$.

Let $A(t) = h(t)S(t)~\forall t\in G$. Fix $s \in G$. We will show that
$L_s A(t) = A(t)$. If this is not true, then there exists a set $E$ of
positive finite measure such that $A(st) \neq A(t)~\forall t\in E$.

The sets $s\Phi^{-1}(E)$ and $\Phi^{-1}(sE)$ may be of $\sigma$-finite
measure. But by choosing a suitable subset $E$ still of positive finite
measure, we can assume that $s\Phi^{-1}(E)$ and $\Phi^{-1}(sE)$ are of
positive finite measure. Having done this, let $F=
s\Phi^{-1}(E)\bigcup\Phi^{-1}(sE)$. Then $\forall t\in E$, $st\in sE
\subseteq \Phi(F)$ and $E\subseteq \Phi(s^{-1}F)$. Now for $t\in E$ and
$x\in X$,
\begin{align*}
L_s(T\chi_{F} x)(t) &= T(\chi_{F} x)(st) \\[.2pc]
&= \chi_{\Phi (F)}(st) A(st)(x)\\[.2pc]
&= A(st)(x).
\end{align*}
Also,
\begin{align*}
T(L_s \chi_{F} x)(t) &= T(\chi_{s^{-1}F} x)(t)\\[.2pc]
&= \chi_{\Phi(s^{-1}F)}(t) A(t)(x)\\[.2pc]
&= A(t)(x).
\end{align*}
Since $L_sT = T L_s$, it follows that $A(st)(x) = A(t)(x)$ for almost
all $t \in E$. By choosing a countable dense set $\{x_n\}_{n=1}^\infty$
in $X$, we conclude that
\begin{equation*}
A(st)(x) = A(t)(x)
\end{equation*}
for almost all $t \in E$ and all $x\in X$. But this is a contradiction.
Hence
\begin{equation*}
A(st) = A(t)
\end{equation*}
for almost all $t \in G$. Therefore for each $x\in X$,
\begin{equation*}
h(t) S(t)(x) = h(st) S(st)(x)
\end{equation*}
for almost all $t \in G$. Since $S(t)$ is an isometry of $X$ onto itself
and $h(t)\geq 0$, we have
\begin{equation*}
h(st) = h(t)
\end{equation*}
for almost all $t \in G$. This implies that $h$ is a constant, say $k$.
It also follows that
\begin{equation*}
S(st) = S(t)
\end{equation*}
for almost all $t \in G$. Hence $S$ is also a constant operator, say
$V$. Therefore, $T$ is an isometric multiplier of $L^p(G,X)$ onto itself
for all $p$, in particular for $p = 1$. Now fix $x \in X$ such that $\|
x\| =1$. Then for $f \in L^{1} (G)$,
\begin{equation*}
L_s T(f x) = L_s k V \Phi(f) x = L_s \Phi(f) k V(x)
\end{equation*}
and
\begin{equation*}
T L_s(f x) = k V \Phi(L_s f) x = \Phi(L_s f) k V(x).
\end{equation*}
Hence $L_s\Phi(f) = \Phi(L_s(f))$. This implies that the map $f
\longrightarrow k \Phi(f)$ is an isometric multiplier of $L^1(G)$ onto
itself. Hence by Wendel's characterization there exists an $s\in G$ and
a scalar $c$ such that $|c|=1$ for which we have
\begin{equation*}
k \Phi(f)(t) = c f (ts).
\end{equation*}
Let $U = k V$. Then $U$ is an isometry of $X$ onto itself such that $T =
U\circ R_s$ and
\begin{equation*}
(TF)(t) = U F(ts)
\end{equation*}
for almost all $t \in G$ and all $F\in L^p(G,X)$. This completes the
proof of the theorem.}\hfill$\Box$
\end{pot}

We shall now show that the condition that $X$ is not an $\ell^p$-direct
sum is a necessary (as well as a sufficient) condition for the
conclusion of Theorem~1 to hold. In fact, we prove the following
theorem.

\begin{theor}[\!] Let $X$ be a separable Banach space which is
$\ell^p$-direct sum of two non-zero subspaces of $X$. Then there exists
an isometric{\rm ,} invertible{\rm ,} left multiplier $T$ of $L^p(G,X)$
which is not of the form $U\circ R_y$ for any isometry $U$ of $X$ and
$y\in G$.
\end{theor}

\begin{proof} Suppose $X = X_1 \oplus_p X_2$. Then
\begin{equation*}
L^p(G,X) = L^p(G,X_1) \oplus_p  L^p(G,X_2).
\end{equation*}
Choose $z\in G$ where $z$ is not the identity element of $G$. Define $T$
by
\begin{equation*}
T(f_1 \oplus f_2) = f_1 \oplus R_z f_2.
\end{equation*}
Then it is easy to verify that $T$ is an isometric, invertible, left
multiplier of $L^p(G,X)$ which is not of the form $U\circ R_y$ for any
isometry $U$ of $X$ and $y\in G$.\hfill$\Box$
\end{proof}

\section{Isometric multipliers of $\pmb{L^1\cap L^p(G,X)}$ and $\pmb{L^1\cap
C_0(G,X)}$}

In this section we assume that $G$ is non-compact and $X$ is not an
$\ell^p$-direct sum of two non-zero subspaces of $X$. We will prove that
if $T$ is an isometric, invertible, left multiplier of $L^1\cap
L^p(G,X)$ or $L^1\cap C_0(G,X)$ then $T$ is of the form $U\circ R_y$ for
some isometry $U$ of $X$ and $y\in G$.

The proof of the following proposition is similar to the proof of
Theorems~3.5.1 and 3.5.2 in \cite{2} and hence omitted.

\setcounter{propo}{2}
\begin{propo}$\left.\right.$\vspace{.5pc}

\noindent Suppose $G$ is non-compact. If $T$ is a left multiplier of
$L^1\cap L^p(G,X)$ or $L^1\cap C_0(G,X)$ then $T$ has a unique extension
to $L^1(G,X)$ as a left multiplier such that $\|Tf\|_1 \leq \|T\|
\|f\|_1${\rm ,} where $\|T\|$ is the norm of $T$ as an operator on
$L^1\cap L^p(G,X)$ or $L^1\cap C_0(G,X)$.
\end{propo}
We now prove the characterization of an isometric, invertible, left
multiplier of $L^1\cap L^p(G,X)$ or $L^1\cap C_0(G,X)$.

\setcounter{theor}{3}
\begin{theor}[\!] Suppose $G$ is non-compact and $X$ is not
$\ell^p$-direct sum of two non-zero subspaces of $X$. If $T$ is an
isometric{\rm ,} invertible{\rm ,} left multiplier of $L^1\cap L^p(G,X)$
or $L^1\cap C_0(G,X)$ then $T$ is of the form $U\circ R_y$ for some
isometry $U$ of $X$ and $y\in G$.
\end{theor}

\begin{proof} Since $T$ and $T^{-1}$ are both isometric multipliers of
$L^1\cap L^p(G,X)$ or $L^1\cap C_0(G,X)$, it follows from Proposition~3
that $T$ extends to $L^1(G,X)$ as an isometric left multiplier.
Therefore by Theorem~1, there exists an isometry of $X$ onto itself and
$y\in G$ such that $T = U\circ R_y$. \hfill $\Box$
\end{proof}

\section{Isometric multipliers of $\pmb{A^p(G,H)}$}

Let $G$ be a locally compact Abelian group and $H$ be a separable
Hilbert space. We define the Fourier transform of $f\in L^1(G,H)$ by
\begin{equation*}
\hat{f}(\gamma) = \int_G \overline{\gamma(x)} f(x) {\rm d}x,
\end{equation*}
where $\gamma\in \Gamma$, the dual group of $G$. Given a Haar measure on
$G$ there exists a unique Haar measure on $\Gamma$ such that the map
$f\rightarrow \hat{f}$ is an isometry of $L^1\cap L^2(G,H)$ into $L^2
(\Gamma,H)$ and extends to an isometry of $L^2 (G,H)$ onto $L^2
(\Gamma,H)$. The Fourier--Plancherel formula $\|\hat{f}\|_2 = \|f\|_2$
holds for $f\in L^2 (G,H)$, see [1].

For $f\in A^p(G,H)$, we define $\|f\| = \|f\|_1 + \|\hat{f}\|_p$. Then
$A^p(G,H)$ is a Banach space. We note that left and right translates
mean the same for Abelian groups. Suppose $G$ is non-compact. We will
prove that if $T$ is an isometric and invertible multiplier of $A^p(G,H)$
then $T = U\circ R_y$, where $U$ is an isometry of $H$ onto itself and
$y\in G$.

The proof of the following Proposition is similar to the argument in the
proof of Theorem~6.3.1 in [2] where it is shown that if $T$ is a multiplier of
$A^p(G)$ then $\|Tf\|_1 \leq \|T\| \|f\|_1$ for $f\in A^p(G)$, where
$\|T\|$ denotes the operator norm of $T$. The necessary modifications
are easy and hence we omit the details.

\setcounter{propo}{4}
\begin{propo}$\left.\right.$\vspace{.5pc}

\noindent Let $G$ be a non-compact locally compact Abelian group and
$1\leq p<\infty$. Suppose $T$ is a multiplier of $A^p(G,H)$ then
$\|Tf\|_1 \leq \|T\| \|f\|_1$ for $f\in A^p(G,H)$.
\end{propo}

We now prove the characterization of isometric multipliers of $A^p(G,H)$.

\setcounter{theor}{5}
\begin{theor}[\!]
Let $G$ be a non-compact locally compact Abelian group and $1\leq
p<\infty$. Suppose $T$ is an isometric multiplier of $A^p(G,H)$. Then
there exists a unique $y\in G$ and an isometry $U$ of $H$ onto itself such
that $T = U\circ R_y$.
\end{theor}

\begin{proof} Let $T$ be an isometric multiplier of $A^p(G,H)$. Then
$T^{-1}$ is also an isometric multiplier and we conclude from
Proposition~5 that $\|Tf\|_1 = \|f\|_1$ for every $f\in A^p(G,H)$. It
follows that $T$ extends to $L^1(G,H)$ as an isometric multiplier of
$L^1(G,H)$. Hence, by Theorem~1, there exists an isometry $U$ of $H$
onto itself and $y\in G$ such that $T = U\circ R_y$.

\hfill$\Box$
\end{proof}

\section{Isometric multipliers of $\pmb{C(G,X)}$}

In this section we describe the isometric, invertible, left multipliers
of $C(G,X)$ where $G$ is a compact group and $X^*$ is strictly convex.
The space $C(G,X)$ consists of all continuous $X$-valued function and is
a Banach space under the supremum norm. The norm of $f\in C(G,X)$ will
be denoted by $\|f\|_\infty$. For the space $X$, we denote the set of
isometries of $X$ onto itself by $I(X)$. The isometries of $C(G,X)$ were
characterized by Lau \cite{3}. He has shown that if $T$ is an isometry
of $C(G,X)$ onto itself, then there exists a homeomorphism $\phi$ of $G$
onto itself and a continuous map $\lambda\hbox{:}\ X\rightarrow I(X)$
(with the strong operator topology) such that
\begin{equation*}
Tf(t) = \lambda(t) f(\phi(t)).
\end{equation*}
Using this characterization of isometries of $C(G,X)$, we prove the
following:

\begin{theorre}{\it
Let $T$ be an isometric{\rm ,} invertible{\rm ,} left multiplier of
$C(G,X)$. Then there exists an isometry $U$ of $X$ onto itself and $y\in G$
such that $T = U\circ R_y$.}
\end{theorre}

\begin{proof} Since $T$ is an isometry of $C(G,X)$, there exists a
continuous map $\lambda\hbox{:}\ X\rightarrow I(X)$ and a homeomorphism
$\phi$ of $G$ onto itself such that
\begin{equation*}
Tf(s) = \lambda(s) f(\phi(s))\quad \forall s\in G.
\end{equation*}
Fix $x\in X$ and let $ f(s) = x~\forall s \in G$. Then
\begin{equation}
T L_t f(s) = \lambda(s) f(t(\phi(s))
\end{equation}
and
\begin{equation}
L_t Tf(s) = \lambda (ts) f(\phi(ts)).
\end{equation}
Since $T L_t = L_t T$, it follows that $\lambda(s)(x) = \lambda(ts)(x)$.
Since $x\in X$ is arbitrary, we conclude that $\lambda(ts) =
\lambda(s)~\forall s, t\in G$. Hence there exists an isometry $U$ of $X$
onto itself such that $\lambda(s) = U~\forall s\in G$. Therefore
\begin{equation*}
Tf(s) = U f(\phi(s))\quad \forall f\in C(G,X).
\end{equation*}
Let $g\in C(G)$ and $x\in X$. Define $f$ by $f(s) = g(s)x~\forall
s\in G$. Then (1) and (2) imply that
\begin{equation*}
g(t\phi(s)) = g(\phi(ts))\quad \forall g \in C(G).
\end{equation*}
Since $C(G)$ separates points, we conclude that $t\phi(s) =
\phi(ts)~\forall s,t \in G$. Let $s$ be the identity element of $G$.
Then $\phi(t) = t\phi(e)$. Let us denote $\phi(e)$ by $y$. Then we have
$Tf(s) = U f(sy)~\forall f\in C(G,X)$ and $s\in G$. Therefore we have $T
= U\circ R_y$.\hfill$\Box$
\end{proof}

\end{document}